\documentclass{amsart}
\usepackage{amssymb}
\usepackage{amsmath}
\usepackage{amsfonts}

\newtheorem{theorem}{Theorem}
\theoremstyle{plain}

\newtheorem{condition}{Condition}

\newtheorem{corollary}{Corollary}

\newtheorem{definition}{Definition}
\newtheorem{example}{Example}

\newtheorem{lemma}{Lemma}

\newtheorem{proposition}{Proposition}
\newtheorem{remark}{Remark}

\numberwithin{equation}{section}
\input{tcilatex}

\begin{document}
\title[Densities for RDEs under H\"{o}rmander's Condition]{Densities for
Rough Differential Equations under H\"{o}rmander's Condition}
\author{Thomas Cass and Peter Friz}

\begin{abstract}
We consider stochastic differential equations $dY=V\left( Y\right) dX$
driven by a multidimensional Gaussian process $X$ in the rough path sense.
Using Malliavin Calculus we show that $Y_{t}$ admits a density for $t\in
(0,T]$ provided (i) the vector fields $V=\left( V_{1},...,V_{d}\right) $
satisfy H\"{o}rmander's condition and (ii)\ the Gaussian driving signal $X$
satisfies certain conditions. Examples of driving signals include
fractional\ Brownian motion with Hurst parameter $H>1/4$, the Brownian
Bridge returning to zero after time $T$ and the Ornstein-Uhlenbeck process.
\end{abstract}

\maketitle

\section{Introduction}

In the theory of stochastic processes, H\"{o}rmander's theorem on
hypoellipticity of degenerate partial differential equations has always been
an important tool to decide whether or not a diffusion process with given
generator admits a density. This dependence on PDE\ theory was removed when
P. Malliavin devised a purely probabilistic approach to H\"{o}rmander's
theorem which is perfectly adapted to prove existence and smoothness of
densities for diffusions given as strong solution to an It\^{o} stochastic
differential equation driven by Brownian motion.

The key ingredients of Malliavin's machinery, better known as \textit{%
Malliavin Calculus} or \textit{stochastic calculus of variations} can be
formulated in the setting of an abstract Wiener space $\left( W,\mathcal{H}%
,\mu \right) $. This concept is standard (e.g. \cite{Ma} or any modern book
on stochastic analysis) as is the notion of \textit{weakly non-degenerate} $%
\mathbb{R}^{e}$-valued functional $\varphi $ which has the desirable
property that the image measure $\varphi _{\ast }\mu $ is absolutely
continuous with respect to Lebesgue measure on $\mathbb{R}^{e}$.
(Functionals which are \textit{non-degenerate }have a smooth density.)
Precise definitions are given later on in the text.

\ 

Given these abstract tools, we turn to the standard Wiener space $C\left( %
\left[ 0,T\right] ,\mathbb{R}^{d}\right) $ equipped with Wiener measure i.e.
the standard model for Brownian motion $B=B\left( \omega \right) $.\ From It%
\^{o}'s theory, we know how to solve the stochastic differential equation%
\begin{equation*}
dY=\sum_{i=1}^{d}V_{i}\left( Y\right) \circ dB^{i}\equiv V\left( Y\right)
\circ dB,\,\,\,Y\left( 0\right) =y_{0}\in \mathbb{R}^{e}.
\end{equation*}%
The It\^{o}-map $B\mapsto Y$ is notorious for its lack of strong regularity
properties. On the positive side, it is smooth in a weak Sobolev type sense
("smooth in Malliavin's sense") and under H\"{o}rmander's condition at $%
y_{0}\in \mathbb{R}^{e}$ 
\begin{equation}
\left( H\right) :\text{ \ Lie}\left[ V_{1},...,V_{d}\right] _{y_{0}}=%
\mathcal{T}_{y_{0}}\mathbb{R}^{e}\cong \mathbb{R}^{e}  \label{CondH}
\end{equation}%
one can show (e.g. \cite{Ma, Sh, Be, Nu06}) that the solution map $B\mapsto
Y_{t}$ is non-degenerate for all $t\in (0,T]$. This line of reasoning
provides a direct probabilistic approach to the study of transition
densities of $Y$ and has found applications from stochastic fluid dynamics
to interest rate theory, e.g. \cite{FT, HM}. The same range of applications%
\footnote{%
For instance, stochastic differential equations driven by fBM have
applications applications to vortex filaments; applications to finance
include (geometric) fractional Brownian motion as paradigm of a
non-semimartingale which admits no arbitrarge under transaction costs. The
reader is refered to the books \cite[Sec. 5.3 and 5.4]{Nu06}, \cite[Sec 8.1.]%
{DeSc06} and the references therein.} nowadays demand stochastic models of
type%
\begin{equation}
dY=V\left( Y\right) dX  \label{GaussRDE}
\end{equation}%
where $X$ is a Gaussian process, such as fractional Brownian motion (short:
fBm). Differential equations of this type have also been used as simple
examples for the study of ergodicity of non-Markovian systemems, \cite{H05}.

In a previous paper \cite{CFV} we linked rough paths and Malliavin calculus
by giving a (simple) proof that existence of a density for solutions of (\ref%
{GaussRDE}) holds true under ellipticity i.e. 
\begin{equation*}
\left( E\right) :\text{ \ Span}\left[ V_{1},...,V_{d}\right] _{y_{0}}=%
\mathcal{T}_{y_{0}}\mathbb{R}^{e}\cong \mathbb{R}^{e}.
\end{equation*}
and generic non-degeneracy conditions on $X$, the differential equation (\ref%
{GaussRDE}) being understood in the rough path sense \cite{L98, LQ02}, a
unified framework which covers \textit{at once} Young and Stratonovich
solutions\ (and goes well beyond). The aim of this paper is to prove the
existence of densities under H\"{o}rmander's condition (H) in the following
form:

\begin{theorem}
\label{Main}Let $\left( X_{t}^{1},\dots ,X_{t}^{d}\right) =(X_{t}:t\in \left[
0,T\right] )$ be a continuous, centered Gaussian process with independent
components $X^{1},\dots ,X^{d}$. Assume $X$ satisfies the conditions listed
in section \ref{Conditions}. (In particular, $X$ is assumed to lift to a
geometric rough path so that (\ref{GaussRDE}) makes sense as random rough
differential equations.) Let $V=(V_{1},...,V_{d})$ a collection of smooth
bounded vector fields on $\mathbb{R}^{e}$ with bounded derivatives which
satisfies H\"{o}rmander's condition (H) at $y_{0}$. Then the random RDE
solution $Y_{t}=Y_{t}\left( \omega \right) \in \mathbb{R}^{e}$ to (\ref%
{GaussRDE}) started at $Y_{0}=y_{0}$ admits a density with respect to
Lebesgue measure on $\mathbb{R}^{e}$ for all times $t\in (0,T].$
\end{theorem}

One should note that $X$, the Gaussian driving signal of (\ref{GaussRDE}),
is fully described by the covariance function of each component and, under
the further assumption of IID\ components, by the covariance of a single
component, i.e. $R\left( s,t\right) =\mathbb{E}\left(
X_{s}^{1}X_{t}^{1}\right) $. In principle all conditions on $X$ are
checkable from the covariance, in practice it is convenient to have
conditions available which involve the \textit{reproducing kernel Hilbert} $%
\sim $ or \textit{Cameron-Martin space} associated to $X$ as well as certain
sample path properties. Leaving these technical details to section \ref%
{Conditions} we emphasize that our conditions are readily checked in many
cases and now give a list of examples to which our theorem applies. It may
be helpful to note that whenever $X$ is a semi-martingale on $\left[ 0,T%
\right] $ then (\ref{GaussRDE}) can be understood as Stratonovich stochastic
differential equation, i.e.%
\begin{equation*}
dY=\sum_{i=1}^{d}V_{i}\left( Y\right) \circ dX^{i}.
\end{equation*}%
In such cases, rough path theory appears as intermediate tool that is
neither needed to understand the assumptions nor the conclusions of theorem %
\ref{Main}. There may be cases when $X$ can be written in terms of Brownian
motion so that ultimately the techniques of \cite{CM, Ta} are applicable.
But in general theorem \ref{Main} covers new grounds.

\begin{example}[Brownian motion]
\label{ExBM} When $R\left( s,t\right) =\min \left( s,t\right) $ the driving
signal $X=X\left( \omega \right) $ is a $d$-dimensional standard Brownian
motion on $\left[ 0,T\right] $. As is well understood \cite{LQ02} one needs
to add L\'{e}vy's area process to obtain a geometric rough path (known in
this context as Brownian rough path or Enhanced Brownian motion). A solution
to (\ref{GaussRDE}) in the rough path sense then precisely solves the
stochastic differential equation in Stratonovich form%
\begin{equation*}
dY=\sum_{i=1}^{d}V_{i}\left( Y\right) \circ dB^{i}.
\end{equation*}%
Subject to H\"{o}rmander's condition (H), Theorem \ref{Main} then shows that 
$Y_{t}=Y_{t}\left( \omega \right) $ has a density for $t>0$ which is of
course well-known.
\end{example}

\begin{example}[Fractional Brownian motion]
\label{ExfBM}When $2R\left( s,t\right) =t^{2H}+s^{2H}-\left\vert
t-s\right\vert ^{2H}$ the corresponding process is known as fractional
Brownian motion $B^{H}=B^{H}\left( \omega \right) $ with Hurst parameter $%
H\in \left( 0,1\right) $. Popularized by \cite{Ma68} it generalizes Brownian
motion which corresponds to $H=1/2$. For fixed $H$ and when no confusion is
possible we shall write $B=(B^{1},\dots ,B^{d})$ for $d$-dimensional fBm.
When $H>1/2$, Kolmogorov's criterion shows that $B$ has nice sample paths
(more precisely, H\"{o}lder continuous sample paths of exponent greater $1/2$%
) which has the great advantage that (\ref{GaussRDE}) can be understood as
integral equation for fixed $\omega $ based on Young integrals (i.e. limits
of Riemann-Stieltjes sums). In this setting of nice sample paths, existence
of a density was established in \cite{NuSa} assuming ellipticity. Using
deterministic estimates for the Jacobian of the flow this density was then
shown to be smooth \cite{NuHu}; building on the same estimates\ the H\"{o}%
rmander case was obtained in \cite{BH}. For $H\leq 1/2$ the situation
appears to be fundamentally different: first, in view of Brownian (and
worse) sample path regularity one needs It\^{o} or rough path ideas to make
sense of (\ref{GaussRDE}) for $H\leq 1/2$. Secondly, the proof of \cite{BH}
does not extend to the rough path setting\footnote{%
The estimates of \cite{NuHu} can be generalized \cite{FVforth} to (sharp)
deterministic estimates on the Jacobian of RDE solutions giving $L^{p}$%
-estimates on the flow of RDEs driven by fBM if and only if $H>1/2$. In
particular, one sees that $L^{p}$-estimates on the flow of Stratonovich SDEs
($H=1/2$) are fundamentally probabilistic i.e. rely on cancellations in
stochastic integration. At present, the question of how to obtain good
integrability when $H<1/2$ is open although one suspects that Gaussian
isoperimetry will ultimately play a role.} and relies somewhat delicately on
specific properties of fractional Brownian motion. In any case, theorem \ref%
{Main} shows that $Y_{t}=Y_{t}\left( \omega \right) $, solution to the RDE
driven by multidimensional fBM with Hurst parameter $H>1/4$ has a density
for all positive times provided the vector fields satisfy H\"{o}rmander's
condition\footnote{%
As is well understood \cite{LQ02}, for $H\leq 1/4$ fractional Brownian
increments decorrelate too slowly for stochastic area to exist and so there
is no meaningful lift of fBM with $H\leq 1/4$ to a geometric rough path.}.
The novelity is of course the degenerate regime $H<1/2$ with sample path
regularity worse than Brownian motion.
\end{example}

\begin{example}[Brownian Bridge]
\label{ExBB}Let $B$ be a $d$-dimensional standard Brownian motion. Define
the Brownian bridge returning to zero at time $T$ by%
\begin{equation*}
X_{t}^{T}:=B_{t}-\frac{t}{T}B_{T}\text{ for }t\in \left[ 0,T\right] .
\end{equation*}%
Equivalently, one can define $X^{T}$ via the covariance%
\begin{equation*}
R^{T}\left( s,t\right) =\min \left( s,t\right) \left( 1-\max \left(
s,t\right) /T\right) .
\end{equation*}%
Clearly $X_{t}^{T}|_{t=T}=0$ and trivially (take $dY=dX$) the conclusion of
Theorem \ref{Main} cannot hold; this behaviour is indeed ruled out by
condition \ref{NonDegGauss}\ in section \ref{Conditions}. On the other hand,
we may consider $X^{T+\varepsilon }$ restricted to $\left[ 0,T\right] $ and
in this case the conditions in section \ref{Conditions} are readily
verified. It is worth remarking that $Z:=X^{T+\varepsilon }$ stopped at time 
$T$ is also a semi-martingale; for instance, by writing $\left(
X_{t}^{T+\varepsilon }:t\leq T\right) $ as strong solution to an It\^{o}
differential equation with (well-behaved) drift (as long as $t\leq T$). The
conclusion of theorem \ref{Main} can then be stated by saying that the
unique Stratonovich solution to $dY=\sum V_{i}\left( Y\right) \circ dZ^{i}$
admits a density for all times $t\in (0,T]$ provided the vector fields
satisfy H\"{o}rmander's condition (H).
\end{example}

\begin{example}[Ornstein-Uhlenbeck]
\label{ExOU}Let $B$ be a standard \thinspace $d$-dimensional Brownian motion
and define the centered Gaussian process $X$ by Wiener-It\^{o} integration,%
\begin{equation*}
X_{t}^{i}=\int_{0}^{t}e^{-\left( t-r\right) }dB_{r}^{i}\text{ with }%
i=1,...,d.
\end{equation*}%
$X$ satisfies the It\^{o} differential equations, $dX_{t}=-X_{t}dt+dB_{t}$
and is also a semi-martingale. The conditions of section \ref{Conditions}
are readily checked (in essence one uses $X_{t}\sim B_{t}$ at $t\rightarrow
0+$ and the absence of Brownian Bridge type degeneracy). The conclusion of
theorem \ref{Main} can then be stated by saying that the unique Stratonovich
solution to $dY=\sum V_{i}\left( Y\right) \circ dX^{i}$ admits a density for
all positive times provided the vector fields satisfy H\"{o}rmander's
condition (H).
\end{example}

Further examples (for instance, "fractional" version of the Brownian Bridge
and Ornstein-Uhlenbeck process) are readily constructed. Generalizing
examples \ref{ExfBM}, \ref{ExOU} one could consider Volterra processes \cite%
{De04}, i.e. Gaussian process with representation $X=\int K\left( \cdot
,s\right) dB_{s}$, and derive sufficient conditions on the kernel $K$ which
imply those of section \ref{Conditions}. Existence of a rough path lift of $X
$ aside, one would need non-degeneracy of $K$ and certain scaling properites
as $t\rightarrow 0+$ but we shall not pursue this here. (In any case, there
are non-Volterra examples, such as the Brownian bridge returning to zero at $%
\left( T+\varepsilon \right) $, to which theorem \ref{Main} applies.)\newline

The proof of theorem \ref{Main} is based on the fact \cite{CFV} that RDE
solutions driven by Gaussian signals are "$\mathcal{H}$-differentiable" i.e.
differentiable in Cameron-Martin directions. Existence of a density is then
reduced to showing that the Malliavin covariance matrix is weakly
non-degenerate. The standard proof of this (e.g. \cite{Ma, Be} or \cite[Sec
2.3.2]{Nu06}) is based on Blumenthal's 0-1 law and the Doob-Meyer
decomposition for semi-martingales. The main difficulty to overcome in the
general Gaussian context of this paper is that the Doob-Meyer decomposition
is not available and we manage to bypass its use by suitable \textit{small
time developments for RDEs}, obtained in \cite{FV07c}; in conjunction with 
\textit{(Stroock-Varadhan type) support description for certain Gaussian
rough paths} (as conjectued by Ledoux et al. \cite{LQZ} and carried out
independently in \cite{FdeP06, FV2}, see also \cite{CFV} and \cite{FV07b}.)

The crucial induction step - which explains the appearance of higher
brackets - requires us to assume a "non-standard" H\"{o}rmander condition
which involves only iterated Lie-brackets contracted against certain tensors
arising from free nilpotent Lie groups. Equivalence to the usual H\"{o}%
rmander condition (H) is then established separately. 

\section{Preliminaries on\ ODE and RDEs}

\subsection{Controlled ordinary differential equations}

Consider the controlled ordinary differential equations, driven by a smooth $%
R^{d}$-valued signal $f=f\left( t\right) $ along sufficiently smooth and
bounded vector fields $V=\left( V_{1},...,V_{d}\right) $,%
\begin{equation}
dy=V(y)df\equiv \sum_{i=1}^{d}V_{i}\left( y\right) f^{\prime }\left(
t\right) dt,\,y\left( t_{0}\right) =y_{0}\in \mathbb{R}^{e}.  \label{cODE}
\end{equation}%
We call $U_{t\leftarrow t_{0}}^{f}\left( y_{0}\right) \equiv y_{t}$ the
associated flow. Let $J$ denote the Jacobian of $U$. It satisfies the ODE
obtain by formal differentiation w.r.t. $y_{0}$. More specifically, 
\begin{equation*}
a\mapsto \left\{ \frac{d}{d\varepsilon }U_{t\leftarrow t_{0}}^{f}\left(
y_{0}+\varepsilon a\right) \right\} _{\varepsilon =0}
\end{equation*}%
is a linear map from $\mathbb{R}^{e}\rightarrow \mathbb{R}^{e}$ and we let $%
J_{t\leftarrow t_{0}}^{f}\left( y_{0}\right) $ denote the corresponding $%
e\times e$ matrix. It is immediate to see that 
\begin{equation*}
\frac{d}{dt}J_{t\leftarrow t_{0}}^{f}\left( y_{0}\right) =\left[ \frac{d}{dt}%
M^{f}\left( U_{t\leftarrow t_{0}}^{f}\left( y_{0}\right) ,t\right) \right]
\cdot J_{t\leftarrow t_{0}}^{f}\left( y_{0}\right)
\end{equation*}%
where $\cdot $ denotes matrix multiplication and%
\begin{equation*}
\frac{d}{dt}M^{f}\left( y,t\right) =\sum_{i=1}^{d}V_{i}^{\prime }\left(
y\right) \frac{d}{dt}f_{t}^{i}.
\end{equation*}%
Note that $J_{t_{2}\leftarrow t_{0}}^{f}=J_{t_{2}\leftarrow t_{1}}^{f}\cdot
J_{t_{1}\leftarrow t_{0}}^{f}.$ We can also consider Gateaux derivates in
the driving signal and define%
\begin{equation*}
D_{h}U_{t\leftarrow 0}^{f}=\left\{ \frac{d}{d\varepsilon }U_{t\leftarrow
0}^{f+\varepsilon h}\right\} _{\varepsilon =0}.
\end{equation*}%
One sees that $D_{h}U_{t\leftarrow 0}^{f}$ satisfies a linear ODE and the
variation of constants formula leads to%
\begin{equation*}
D_{h}U_{t\leftarrow 0}^{f}\left( y_{0}\right)
=\int_{0}^{t}\sum_{i=1}^{d}J_{t\leftarrow s}^{f}\left( V_{i}\left(
U_{s\leftarrow 0}^{f}\right) \right) dh_{s}^{i}.
\end{equation*}%
Finally, given a smooth vector field $W$ a straight-forward computation
gives 
\begin{equation}
dJ_{0\leftarrow t}^{f}\left( W\left( U_{t\leftarrow 0}^{f}\right) \right)
=J_{0\leftarrow t}^{f}\left( \left[ V_{i},W\right] \left( U_{t\leftarrow
0}^{f}\right) \right) df_{t}^{i}.  \label{LieIdentity}
\end{equation}

\subsection{Rough differential equations}

Following \cite{L98, LQ02, FV07} a geometric $p$-rough path $\mathbf{x}$
over $\mathbb{R}^{d}$ is a continuous path on $\left[ 0,T\right] $ with
values in $G^{\left[ p\right] }\left( \mathbb{R}^{d}\right) $, the step-$%
\left[ p\right] $ nilpotent group over $\mathbb{R}^{d}$, and of finite $p$%
-variation relative to the \cite{FV07} Carnot-Caratheodory metric $d$ on $G^{%
\left[ p\right] }\left( \mathbb{R}^{d}\right) $, i.e.%
\begin{equation*}
\sup_{n\in \mathbb{N}}\sup_{0<t_{1}<...<t_{n}<T}\sum_{i}d\left( \mathbf{x}%
_{t_{i}},\mathbf{x}_{t_{i+1}}\right) ^{p}<\infty .
\end{equation*}%
As in \cite{chen, L98} we view $G^{\left[ p\right] }\left( \mathbb{R}%
^{d}\right) $ as embedded in its enveloping tensor algebra i.e.%
\begin{equation*}
G^{\left[ p\right] }\left( \mathbb{R}^{d}\right) \subset T^{\left[ p\right]
}\left( \mathbb{R}^{d}\right) :=\oplus _{k=0,...,\left[ p\right] }\left( 
\mathbb{R}^{d}\right) ^{\otimes k}\text{.}
\end{equation*}
One can then think of $\mathbf{x}$ as a path $x:\left[ 0,T\right]
\rightarrow \mathbb{R}^{d}$ enhanced with its iterated integrals although
the later need not make classical sense\footnote{%
In fact, $G^{N}\left( \mathbb{R}^{d}\right) $ can realized as all points in
the tensor algebra which arise from computing iterated integrals up to order 
$N$ of smooth paths over a fixed time interval. The group product then
corresponds to the concatenation of paths, the inverse corresponds to
running a path backwards in time etc.}. The canonical projection to $\left( 
\mathbb{R}^{d}\right) ^{\otimes k}$ is denoted $\pi _{k}\left( \mathbf{x}%
\right) $ or $\mathbf{x}^{k}$. Lyons' theory of rough paths then gives
deterministic meaning to the rough differential equation (short: RDE)%
\begin{equation}
dy=V\left( y\right) d\mathbf{x}\text{.}  \label{detRDE}
\end{equation}%
(One can think of RDE solutions as limit points of corresponding ODEs of
form (\ref{cODE}) in which the smooth driving signals \textit{plus their
iterated integrals up to order} $\left[ p\right] $ converge to $\mathbf{x}$
in suitable $p$-variation distance.) The motivating example, e.g. \cite{L98,
LQ02}, is that \textit{almost every}\ continuous joint realization of
Brownian motion and L\'{e}vy's area process (equivalently: iterated
Stratonovich integrals) gives rise to a geometric $p$-rough path for $p>2$,
known as Brownian rough path or Enhanced Brownian motion (cf. example \ref%
{ExBM}) which provides in particular a robust path-by-path view of
Stratonovich SDEs.

Back to the deterministic RDE (\ref{detRDE}) and assuming smoothness of the
vector fields $V=\left( V_{1},...,V_{d}\right) $, the solution induces a
flow $y_{0}\mapsto U_{t\leftarrow t_{0}}^{\mathbf{x}}\left( y_{0}\right) $.
Following \cite{LQ97, LQ98}, the Jacobian $J_{t\leftarrow t_{0}}^{\mathbf{x}%
} $ of the flow exists and satisfies a linear RDE, as does the directional
derivative%
\begin{equation*}
D_{h}U_{t\leftarrow 0}^{\mathbf{X}}=\left\{ \frac{d}{d\varepsilon }%
U_{t\leftarrow 0}^{T_{\varepsilon h}\mathbf{x}}\right\} _{\varepsilon =0}
\end{equation*}%
for a smooth path $h$. If $\mathbf{x}$ arises from a smooth path $x$
together with its iterated integrals the \textit{translated rough path} $%
T_{h}\mathbf{x}$ (cf. \cite{LQ97, LQ02}) is nothing but $x+h$ together with
its iterated integrals. In the general case, we assume $h\in C^{q\text{-var}%
} $ with $1/p+1/q>1$, the translation $T_{h}\mathbf{x}$ can be written in
terms of $\mathbf{x}$ and cross-integrals between $\pi _{1}\left( \mathbf{x}%
_{0,\cdot }\right) =:x$ and the perturbation $h$. (These integrals are
well-defined Young-integrals.)

\begin{proposition}
\label{VarOfConst}Let $\mathbf{X}$ be a geometric $p$-rough paths over $%
\mathbb{R}^{d}$ and $h\in C^{q\text{-var}}\left( \left[ 0,T\right] ,\mathbb{R%
}^{d}\right) $ such that $1/p+1/q>1$. Then%
\begin{equation*}
D_{h}U_{t\leftarrow 0}^{\mathbf{X}}\left( y_{0}\right)
=\int_{0}^{t}\sum_{i}J_{t\leftarrow s}^{\mathbf{X}}\left( V_{i}\left(
U_{s\leftarrow 0}^{\mathbf{X}}\right) \right) dh_{s}^{i}
\end{equation*}%
where the right hand side is well-defined as Young intergral.
\end{proposition}

\begin{proof}
$J_{t\leftarrow 0}^{\mathbf{X}},D_{h}U_{t\leftarrow 0}^{\mathbf{X}}$ satisfy
(at least jointly with $U_{t\leftarrow 0}^{\mathbf{X}}$) RDEs driven by $%
\mathbf{X}$ which allows, in essence, to use Lyons' limit theorem; this is
discussed in detail in \cite{LQ97, LQ98}. A little care is needed since the
resulting vector fields are not bounded anymore. However we can rule out
explosion and then localize the problem: the needed remark is that $%
J_{t\leftarrow 0}^{\mathbf{X}}$ also satisfy a linear RDE of form 
\begin{equation*}
dJ_{t\leftarrow 0}^{\mathbf{X}}=dM^{\mathbf{X}}\left( U_{t\leftarrow 0}^{%
\mathbf{X}}\left( y_{0}\right) ,t\right) \cdot J_{t\leftarrow 0}^{\mathbf{X}%
}\left( y_{0}\right) 
\end{equation*}%
and explosion can be ruled out by direct iterative expansion and estimates
of the Einstein sum as in \cite{L98}.
\end{proof}

\section{RDEs driven by Gaussian signals\label{RDEsDrivenByX}}

We consider a continuous, centered Gaussian process $X=\left(
X^{1},...,X^{d}\right) $ with independent components started at zero. This
gives rise to an abstract Wiener space $\left( W,\mathcal{H},\mu \right) $
where $W=\mathcal{\bar{H}}\subset C_{0}\left( \left[ 0,T\right] ,\mathbb{R}%
^{d}\right) $. Note that $\mathcal{H}=\oplus _{i=1}^{d}\mathcal{H}^{\left(
i\right) }$ and recall that element of $\mathcal{H}$ are of form $h_{t}=%
\mathbb{E}\left( X_{t}\xi \left( h\right) \right) $ where $\xi \left(
h\right) $ is a Gaussian random variable. The ("reproducing kernel")
Hilbert-structure on $\mathcal{H}$ is given by $\,\left\langle h,h^{\prime
}\right\rangle _{\mathcal{H}}:=\mathbb{E}\left( \xi \left( h\right) \xi
\left( h^{\prime }\right) \right) $.

Existence of a Gaussian geometric $p$-rough path above $X$ is tantamount to
the existence of certain L\'{e}vy area integrals. The case of fractional
Brownian motion is well understood and several construction have been
carried out \cite{CQ02, LQ02, FdeP06, Pi07}. In particular, one requires $%
H>1/4$ for the existence of stochastic areas (which can be defined as $%
L^{2}\left( \mathbb{P}\right) $-limits as in It\^{o}'s theory). Resultingly,
one has to deal with geometric $p$-rough paths for $p<4$.\ (When $p<2$ there
is enough sample path regularity to use Young integration and we avoid
speaking of rough paths.)

\begin{condition}
\label{CMqvarCond}Assume $X$ lifts to a (random) geometric $p$-rough path $%
\mathbf{X}$ and $\exists q:1/p+1/q>1$ such that \ 
\begin{equation*}
\mathcal{H}\hookrightarrow C^{q\text{-var}}\left( \left[ 0,T\right] ,\mathbb{%
R}^{d}\right) .
\end{equation*}
\end{condition}

The example to have in mind is Brownian motion for which the above condition
is satisfied with $p=2+\varepsilon $ and $q=1$. (We shall say more about
other Gaussian examples in section \ref{Conditions}.)

If $\mathbf{X}=\mathbf{B}^{H}$ denotes the geometric $p$-rough path, $p\in
\left( 1/H,\left[ 1/H\right] +1\right) $, associated to fractional Brownian
motion then it satisfies a \textit{Stroock-Varadhan support description in
rough path topology}. This was first conjectued by Ledoux et al. \cite{LQZ}
(who obtained it for the Brownian rough path) and carried out independently
in \cite{FdeP06, FV2} for $H>1/3$. The case of $H>1/4$ is more difficult and
discussed in some detail in \cite{CFV}. A proof in the generic context of
Gaussian rough paths (covering fBM with $H>1/4$ as special case) is given in
the forthcoming paper \cite{FV07b}. The statement is%
\begin{equation}
\text{supp}\left( \mathbb{P}_{\ast }\mathbf{X}\right) =\overline{\left\{ S_{%
\left[ p\right] }\left( \mathcal{H}\right) \right\} }  \label{SupportEGP}
\end{equation}%
where support and closure are relative to the homogenous $p$-variation
topology for geometric $p$-rough paths. We recall that $S_{\left[ p\right] }$%
, for $\left[ p\right] =2,3$ given by 
\begin{eqnarray*}
S_{2} &:&h\mapsto 1+\int_{0}^{t}dh+\int_{0}^{t}\int_{0}^{s}dh\otimes dh \\
S_{3} &:&h\mapsto 1+\int_{0}^{t}dh+\int_{0}^{t}\int_{0}^{s}dh\otimes
dh+\int_{0}^{t}\int_{0}^{s}\int_{0}^{r}dh\otimes dh\otimes dh
\end{eqnarray*}%
lifts $\mathbb{R}^{d}$-valued paths canonically to $G^{\left[ p\right]
}\left( \mathbb{R}^{d}\right) $-valued paths. In \cite{FV07b} it is seen that%
\textbf{\ }$\mathbf{X}$ exists provided the covariance has finite $\rho $%
-variation with \thinspace $\rho <2$ and it is also established that $%
\mathcal{H}\hookrightarrow C^{\rho \text{-var}}$ which guarantees that $S_{%
\left[ p\right] }\left( \mathcal{H}\right) $ is well-defined via Young
integration. Such support description will be important in checking
condition \ref{CondLikeFBMatZero}, section \ref{Conditions}. 

\begin{definition}
\cite{Sh, Nu06, Ma}Given an abstract Wiener space $\left( W,\mathcal{H},\mu
\right) $, a r.v. $F:W\rightarrow \mathbb{R}$ is $\mathcal{H}$
differentiable at $\omega \in W$ iff exists $DF\left( \omega \right) \in 
\mathcal{H}^{\ast }$such that%
\begin{equation*}
\forall h\in \mathcal{H}:\left\{ \frac{d}{d\varepsilon }F\left( \omega
+\varepsilon h\right) \right\} _{\varepsilon =0}=\left\langle DF\left(
\omega \right) ,h\right\rangle _{\mathcal{H}}.
\end{equation*}%
A vector-valued r.v. $F=\left( F^{1},...,F^{e}\right) :W\rightarrow \mathbb{R%
}^{e}$ is $\mathcal{H}$ differentiable iff each $F^{i}$ is $\mathcal{H}$
differentiable. In this case, $DF\left( \omega \right) =\left( DF^{1}\left(
\omega \right) ,...,DF^{e}\left( \omega \right) \right) $ is a linear
bounded map from $\mathcal{H}\rightarrow \mathbb{R}^{e}$. One then defines
the Malliavin covariance matrix as the random matrix%
\begin{equation*}
\sigma \left( \omega \right) :=\left( \left\langle
DF^{i},DF^{j}\right\rangle _{\mathcal{H}}\right) _{i,j=1,...,e}\in \mathbb{R}%
^{e\times e}.
\end{equation*}%
We call $F$ weakly non-degenerate if $\det \left( \sigma \right) \neq 0$
almost surely.
\end{definition}

\begin{proposition}
\label{GateauxFrechet}Assume condition (\ref{CMqvarCond}). Then, for fixed $%
t\geq 0$, the $\mathbb{R}^{e}$-valued random variable%
\begin{equation*}
\omega \mapsto U_{t\leftarrow 0}^{\mathbf{X}\left( \omega \right) }\left(
y_{0}\right)
\end{equation*}%
is almost surely $\mathcal{H}$-differentiable.
\end{proposition}

\begin{proof}
By assumption $1/p+1/q>1$. We may assume that $\mathbf{X}\left( \omega
\right) $ has been defined so that $\mathbf{X}\left( \omega \right) $ is a
geometric $p$-rough path for every $\omega \in W$. We also know that $%
\mathbf{X}\left( \omega +h\right) \equiv T_{h}\mathbf{X}$ almost surely and
so%
\begin{equation*}
\mathbb{P}\left[ \mathbf{X}\left( \omega +\varepsilon _{n}h\right) \equiv
T_{\varepsilon h}\mathbf{X}\text{ for any countable family }\left(
\varepsilon _{n}\right) \right] =1
\end{equation*}%
We fix $\omega $ in the above set of full measure. For fixed $t$, define 
\begin{equation*}
Z_{i,s}\equiv J_{t\leftarrow s}^{\mathbf{X}}\left( V_{i}\left(
U_{s\leftarrow 0}^{\mathbf{X}}\right) \right) .
\end{equation*}%
Noting that $s\mapsto Z_{i,s}$ is in $C^{p\text{-var}}$ we have, with
implicit summation over $i$, 
\begin{eqnarray*}
\left\vert D_{h}U_{t\leftarrow 0}^{\mathbf{X}}\left( y_{0}\right)
\right\vert &=&\left\vert \int_{0}^{\tau }J_{t\leftarrow s}^{\mathbf{X}%
}\left( V_{i}\left( U_{s\leftarrow 0}^{\mathbf{X}}\right) \right)
dh_{\lambda }^{i}\right\vert \\
&=&\left\vert \int_{0}^{\tau }Z_{i}dh^{i}\right\vert \\
&\leq &c\left( \left\vert Z\right\vert _{p-var}+\left\vert Z\left( 0\right)
\right\vert \right) \times \left\vert h\right\vert _{\rho -var} \\
&\leq &c\left( \left\vert Z\right\vert _{p-var}+\left\vert Z\left( 0\right)
\right\vert \right) \times \left\vert h\right\vert _{\mathcal{H}}
\end{eqnarray*}%
(We used Young's inequality.) The linear map $DU_{t\leftarrow 0}^{\mathbf{X}%
}\left( y_{0}\right) :h\mapsto D_{h}U_{t\leftarrow 0}^{\mathbf{X}}\left(
y_{0}\right) \in \mathbb{R}^{e}$ is bounded and each component is an element
of $\mathcal{H}^{\ast }$, hence%
\begin{equation*}
h\mapsto \left\{ \frac{d}{d\varepsilon }U_{t\leftarrow 0}^{T_{\varepsilon h}%
\mathbf{X}\left( \omega \right) }\left( y_{0}\right) \right\} _{\varepsilon
=0}=\left\langle DU_{t\leftarrow 0}^{\mathbf{X}}\left( y_{0}\right)
,h\right\rangle _{\mathcal{H}}.
\end{equation*}%
Noting that the derivative at $\varepsilon =0$ exists, by definition, if the
difference quotients converge as $\varepsilon \downarrow 0$ and this holds
iff convergence to the same limit takes place along any sequence $%
\varepsilon _{n}\downarrow 0$.\ It follows that, for almost every $\omega $,%
\begin{equation*}
h\mapsto \left\{ \frac{d}{d\varepsilon }U_{t\leftarrow 0}^{\mathbf{X}\left(
\omega +\varepsilon h\right) }\left( y_{0}\right) \right\} _{\varepsilon
=0}=\left\langle DU_{t\leftarrow 0}^{\mathbf{X}}\left( y_{0}\right)
,h\right\rangle _{\mathcal{H}}
\end{equation*}%
and so the random variable $U_{t\leftarrow 0}^{\mathbf{X}}\left(
y_{0}\right) $ is indeed a.s. $\mathcal{H}$-differentiable.
\end{proof}

\section{Conditions on Driving Process\label{Conditions}}

We now give a complete list of assumptions on the ($d$-dimensional) Gaussian
driving signal $\left( X_{t}:t\in \left[ 0,T\right] \right) $. The first
condition was already needed in the previous section to show $\mathcal{H}$%
-differentiability of RDE solutions driven by $X$; we repeat it for
completeness and to give some additional examples.

\begin{condition}
Assume $X$ lifts to a (random) geometric $p$-rough path $\mathbf{X}$ and $%
\exists q:1/p+1/q>1$ such that \ 
\begin{equation}
\mathcal{H}\hookrightarrow C^{q\text{-var}}\left( \left[ 0,T\right] ,\mathbb{%
R}^{d}\right) .  \label{Cond2_Hembds}
\end{equation}
\end{condition}

In the Brownian motion case this holds, as already remarked earlier, with $%
p=2+\varepsilon $ and $q=1$. The same is true for the Brownian bridge and
the Ornstein-Uhlenbeck examples discussed in the introduction; although
case-by-case verifications are not difficult, there is general criterion on
the covariance which implies (\ref{Cond2_Hembds}), see \cite[Prop 16 applied
with $\rho =1$]{FV07}, which also covers fBM. Let us give a direct argument
for case of fBM which covers any Hurst parameter $H>1/4$. Writing $\mathcal{H%
}^{H}$ for the Cameron-Martin space of fBM, the variation embedding in \cite%
{FV3} gives%
\begin{equation*}
\mathcal{H}^{H}\hookrightarrow C^{q\text{-var}}\text{ for any }q>\left(
H+1/2\right) ^{-1}.
\end{equation*}%
At the same time \cite{CQ02, LQ02, FdeP06, Pi07} fBM lifts to a geometric $p$%
-rough path for $p>1/H$. By choosing $p,q$ small enough $1/p+1/q$ can be
made arbitrarily close to $H+\left( H+1/2\right) =2H+1/2>1$ and so (\ref%
{Cond2_Hembds}) holds indeed for fBM with Hurst parameter $H>1/4$.

\begin{condition}
\label{NonDegGauss}Fix $T>0$. We assume non-degeneracy on $\left[ 0,T\right] 
$ in the sense that for any smooth $f=\left( f_{1},...,f_{d}\right) :\left[
0,T\right] \rightarrow \mathbb{R}^{d}$ we have%
\begin{equation*}
\left( \int_{0}^{T}fdh\equiv \sum_{j=1}^{d}\int_{0}^{T}f_{j}dh^{j}=0\forall
h\in \mathcal{H}\right) \implies f\equiv 0.
\end{equation*}
\end{condition}

Again, fBM satisfies the following non-degeneracy condition simply because $%
C_{0}^{1}\left( \left[ 0,T\right] ,\mathbb{R}^{d}\right) \subset \mathcal{H}%
^{H}$, cf. \cite{FV2}. A Brownian bridge which returns to zero at time $T$ \
is ruled out, while a Brownian bridge which returns to zero after time $T$
is allowed. Checking condition \ref{NonDegGauss} for the Ornstein-Uhlenbeck
example in the introduction is left as exercise for the reader. The
following lemma taken from \cite{CFV} contains a few ramifications
concerning condition \ref{NonDegGauss}; since $\mathcal{H}=\oplus _{i=1}^{d}%
\mathcal{H}^{\left( i\right) }$ there is no loss in generality in assuming $%
d=1$.

\begin{lemma}
Assume%
\begin{equation*}
1/p+1/q>1.
\end{equation*}%
(i) The requirement that $f$ is smooth above can be relaxed to $f\in C^{p%
\text{-var}}$.\newline
(ii) The requirement that $\int fdh=0\forall h\in \mathcal{H}$ can be
relaxed to the the quantifier "for all $h$ in some orthonormal basis of $%
\mathcal{H}$".\newline
(iii) The non-degeneracy condition \ref{NonDegGauss} is equivalent to saying
that for all smooth $f\neq 0$, the zero-mean Gaussian random variable $%
\int_{0}^{T}fdX$ (which exists as Young integral or via
integration-by-parts) has positive definite variance.\newline
(iv) The non-degeneracy condition \ref{NonDegGauss} is equivalent to saying
that for all times $0<t_{1}<...<t_{n}<T$ the covariance matrix of $\left(
X_{t_{1}},...,X_{t_{n}}\right) $, that is,%
\begin{equation*}
\left( R\left( t_{i},t_{j}\right) \right) _{i,j=1,...,d}
\end{equation*}%
is (strictly) positive definite. \newline
(v) Non-degeneracy on $\left[ 0,T\right] $ implies non-degeneracy on $\left[
0,t\right] $ for any $t\in (0,T]$.
\end{lemma}

\begin{condition}
\label{ZeroOne}\textbf{"0-1 law":} The germ $\sigma $-algebra $\cap
_{t>0}\sigma \left( X_{s}:s\in \left[ 0,t\right] \right) $ contains only
events of probability zero or one.
\end{condition}

When $X$ is Brownian motion, this is the well-known Blumenthal zero-one law.
More generally, it holds whenever $X$ is an adapted functional of Brownian
motion, including all examples (such as fBM) in which $X$ has a Volterra
presentation \cite{De04}%
\begin{equation*}
X_{t}=\int_{0}^{t}K\left( t,s\right) dB_{s}.
\end{equation*}%
(Nothing is assumed on $K$ other than having the above Wiener-It\^{o}
integral well-defined.) The 0-1 law also holds when $X$ is the strong
solution of an SDE driven by Brownian motion; this includes the
Ornstein-Uhlenbeck - and Browian bridge examples. An example where the 0-1
law fails is given by the \textit{random-ray} $X:t\mapsto tB_{T}\left(
\omega \right) $ in which case the germ-event $\left\{ \omega :dX_{t}\left(
\omega \right) /dt|_{t=0+}\geq 0\right\} $ has probability $1/2$. (In fact,
sample path differentiability at $0+$ implies non-triviality of the germ $%
\sigma $-algebra see \cite{Du} and references therein). We observe that the
random ray example is (a) already ruled out by condition \ref{NonDegGauss}
and (b) should be ruled out anyway since it does not trigger to the bracket
phenomenon needed for a H\"{o}rmander statement.

The next condition expresses some sort of scaled support statement at $t=0+$
and is precisely what is needed in the last part (Step 4) in the proof of
the Theorem \ref{Main} below. We give examples and easier-to-check
conditions below.\ To state it, we recall \cite[Thm 2.2.1]{L98} that a
geometric $p$-rough path $\mathbf{x}$ lifts uniquely lifts uniquely and
continuously (with respect to homogenous $p$-variation distances) to a path
in the free step-$N$ nilpotent group\footnote{%
The $0$ in $C_{0}^{p\text{-var}}$ indicates that $\mathbf{X}_{0}$ is started
at the unit element in the group.}, say%
\begin{equation*}
S_{N}\left( \mathbf{x}\right) \in C_{0}^{p\text{-var}}\left( \left[ 0,T%
\right] ,G^{N}\left( \mathbb{R}^{d}\right) \right) \text{ for }N\geq \left[ p%
\right] .
\end{equation*}%
We also recall that $G^{N}\left( \mathbb{R}^{d}\right) $ carries a dilation
operator $\delta $ which generalizes scalar multiplication on $\mathbb{R}^{d}
$.

\begin{condition}
\label{CondLikeFBMatZero}Assume there exists $H\in \left( 0,1\right) $ such
that for all fixed $N\geq \left[ p\right] $, writing $\mathbf{\tilde{X}}%
=S_{N}\left( \mathbf{X}\right) $, all $g\in $ $G^{N}\left( \mathbb{R}%
^{d}\right) $ and for all $\varepsilon >0$,%
\begin{equation*}
\underset{n\rightarrow \infty }{\lim \inf }\mathbb{P}\left( d\left( \delta
_{n^{H}}\mathbf{\tilde{X}}_{1/n},g\right) <\varepsilon \right) >0.
\end{equation*}
\end{condition}

\begin{proposition}
Let $B$ denote $d$-dimensional fractional Brownian motion with fixed Hurst
parameter $H\in \left( 1/4,1\right) $ and consider the lift to a (random)
geometric $p$-rough path, denoted by $\mathbf{X}=\mathbf{B}$, with $p<4$.
Then it satisfies condition \ref{CondLikeFBMatZero}.
\end{proposition}

\begin{remark}
Brownian motion is covered with $H=1/2$.
\end{remark}

\begin{proof}
Write $\mathbf{\tilde{B}}=S_{N}\left( \mathbf{B}\right) $. From section \ref%
{RDEsDrivenByX}, and the references therein, the support of the law of $%
\mathbf{B}$ w.r.t. homogeneous $p$-variation distance is $C_{0}^{0,p\text{%
-var}}\left( \left[ 0,T\right] ,G^{\left[ p\right] }\left( \mathbb{R}%
^{d}\right) \right) $, that is, the closure of lifted smooth path started at 
$0$ with respect to homogeneous $p$-variation distance \cite{L98, FV07}. By
continuity of $S_{N}$ \cite[Thm 2.2.1]{L98} followed by evaluation of the
path at time $1$ it follows that the support of the law of $\mathbf{\tilde{B}%
}_{1}$ is full, that is, equal to $G^{N}\left( \mathbb{R}^{d}\right) $. On
the other hand, fractional scaling $\left( n^{H}B_{t/n}:t\geq 0\right) 
\overset{D}{=}\left( B_{t}:t\geq 0\right) $ implies $\delta _{n^{H}}\mathbf{%
\tilde{B}}_{1/n}\overset{D}{=}\mathbf{\tilde{B}}_{1}$ and so, thanks to full
support of $\mathbf{\tilde{B}}_{1}$ $,$ 
\begin{equation*}
\underset{n\rightarrow \infty }{\lim \inf }\mathbb{P}\left( d\left( \delta
_{n^{H}}\mathbf{\tilde{B}}_{1/n},g\right) <\varepsilon \right) =\mathbb{P}%
\left( d\left( \mathbf{\tilde{B}}_{1},g\right) <\varepsilon \right) >0.
\end{equation*}
\end{proof}

Although scaling was important in the previous proof, it is only used at
times near $0+$. One thus suspects that every other Gaussian signal $X$
which scales similarly (on the level of $N^{th}$ iterated integrals!)\ also
satisfies condition \ref{CondLikeFBMatZero}. To make this precise we need

\begin{theorem}[\protect\cite{FV07} ]
\label{ThmCiteFV07}Let $\left( X,Y\right) =\left( X^{1},Y^{1},\ldots
,X^{d},Y^{d}\right) $ be a centered continuous Gaussian process on $\left[
0,1\right] $ such that $\left( X^{i},Y^{i}\right) $ are independent for $%
i=1,...,d$. Let $\rho \in \lbrack 1,2)$ and assume the covariance of $\left(
X,Y\right) $, as function on $\left[ 0,1\right] ^{2}$, is of finite $\rho $%
-variation (in 2D sense\footnote{%
Given a function $f$ from $\left[ 0,1\right] ^{2}$ into some normed space,
its variation (in the 2D sense!) is an immediate generalization of the
standard definition but based on "increments" of form%
\begin{equation*}
\Delta _{\lbrack a,b)\times \lbrack c,d)}=f\left( b,d\right) +f\left(
a,c\right) -f(a,d)-f(b,c).
\end{equation*}%
}). Then, for every $p>2\rho $, $X$ and $Y$ can be lifted to geometric $p$%
-rough paths denoted $\mathbf{X}$ and $\mathbf{Y}$. Moreover, there exist a
constant $C$ depending only on $p,\rho $, the covariance of $\left(
X,Y\right) $ so that for all $q\in \lbrack 1,\infty ),$%
\begin{equation*}
\left\vert d_{p\text{-var}}\left( \mathbf{X},\mathbf{Y}\right) \right\vert
_{L^{q}\left( \mathbb{P}\right) }\leq C\sqrt{q}\left\vert R_{X-Y}\right\vert
_{\infty ;\left[ 0,1\right] ^{2}}^{\theta }.
\end{equation*}
\end{theorem}

(Note that $R_{X-Y}\left( s,t\right) $ is a diagonal matrix with entries
depending on $s,t$.)

\begin{corollary}
\label{FV07cor}Let $\left( X,B\right) $ satisfy the conditions of the
previous theorem and assume that $B$ is a ($d$-dimensional) fractional
Brownian motion with fixed Hurst parameter $H\in \left( 1/4,1\right) $.
Assume in addition that%
\begin{equation}
n^{2H}\left\vert R_{X-B}\right\vert _{\infty ;\left[ 0,1/n\right]
^{2}}\rightarrow 0.  \label{CondOfStableCorollary}
\end{equation}%
Then condition \ref{CondLikeFBMatZero} holds.
\end{corollary}

\begin{proof}
With focus on one diagonal entry and with mild abuse of notation (writing $%
X,B$ instead of $X^{i},B^{i}...$)%
\begin{eqnarray*}
&&n^{2H}\left\vert R_{X-B}\right\vert _{\infty ;\left[ 0,1/n\right] ^{2}} \\
&=&\sup_{s,t\in \left[ 0,1\right] }\mathbb{E[}n^{H}\left(
X_{s/n}-B_{s/n}\right) n^{H}\left( X_{t/n}-B_{t/n}\right) ]
\end{eqnarray*}%
which can be rewritten in terms of the rescaled process $X^{\left( n\right)
}=n^{H}X_{\cdot /n}$, and similarly for $B$, as%
\begin{equation*}
\sup_{s,t\in \left[ 0,1\right] }\mathbb{E}\left[ \left( X_{s}^{\left(
n\right) }-B_{s}^{\left( n\right) }\right) \left( X_{t}^{\left( n\right)
}-B_{t}^{\left( n\right) }\right) \right] =\left\vert R_{X^{\left( n\right)
}-B^{\left( n\right) }}\right\vert _{\infty ;\left[ 0,1\right] ^{2}}.
\end{equation*}%
By assumption and the previous theorem, this entails that%
\begin{equation*}
d_{p\text{-var}}\left( \mathbf{X}^{\left( n\right) },\mathbf{B}^{\left(
n\right) }\right) \rightarrow 0\text{ in probability.}
\end{equation*}%
By continuity of $S_{N}$, still writing $\mathbf{\tilde{X}}^{\left( n\right)
}=S_{N}\left( \mathbf{X}^{\left( n\right) }\right) $ for fixed $N$, and
similarly for $\mathbf{B}^{\left( n\right) }$, we have%
\begin{equation*}
d\left( \mathbf{\tilde{X}}_{1}^{\left( n\right) },\mathbf{\tilde{B}}%
_{1}^{\left( n\right) }\right) \leq d_{p\text{-var;}\left[ 0,1\right]
}\left( \mathbf{\tilde{X}}^{\left( n\right) },\mathbf{\tilde{B}}^{\left(
n\right) }\right) \rightarrow 0\text{ in probability.}
\end{equation*}%
But then%
\begin{eqnarray*}
&&\mathbb{P}\left( d\left( \delta _{n^{H}}\mathbf{\tilde{X}}_{1/n},g\right)
<\varepsilon \right) \\
&=&\mathbb{P}\left( d\left( \mathbf{\tilde{X}}_{1}^{\left( n\right)
},g\right) <\varepsilon \right) \\
&\geq &\mathbb{P}\left( d\left( \mathbf{\tilde{X}}_{1}^{\left( n\right) },%
\mathbf{\tilde{B}}_{1}^{\left( n\right) }\right) +d\left( \mathbf{\tilde{B}}%
_{1}^{\left( n\right) },g\right) <\varepsilon \right) \\
&\geq &\mathbb{P}\left( d\left( \mathbf{\tilde{B}}_{1}^{\left( n\right)
},g\right) <\varepsilon /2\right) \\
&&-\mathbb{P}\left( d\left( \mathbf{\tilde{X}}_{1}^{\left( n\right) },%
\mathbf{\tilde{B}}_{1}^{\left( n\right) }\right) >\varepsilon /2\right)
\end{eqnarray*}%
and so%
\begin{equation*}
\underset{n\rightarrow \infty }{\lim \inf }\mathbb{P}\left( d\left( \delta
_{n^{H}}\mathbf{\tilde{X}}_{1/n},g\right) <\varepsilon \right) \geq \,%
\underset{n\rightarrow \infty }{\lim \inf }\mathbb{P}\left( d\left( \mathbf{%
\tilde{B}}_{1}^{\left( n\right) },g\right) <\varepsilon /2\right)
\end{equation*}%
and this is positive by the example in which we discussed the case of $%
\mathbf{B}$ resp. $\mathbf{\tilde{B}}$. The proof is finished.
\end{proof}

\begin{example}[Ornstein Uhlenbeck]
Given a Brownian motion $B$, we consider the Ornstein Uhlenbeck (short:OU)
process given by the It\^{o} integral%
\begin{equation*}
X_{t}=\int_{0}^{t}e^{-\left( t-r\right) }dB_{r}.
\end{equation*}%
If $B$ is $d$-dimensional this yields, component-wise, the $d$-dimensional
OU process. It is readily checked that $\left( X,B\right) $ satisifes the
assumptions of theorem \ref{ThmCiteFV07}. (In fact, one sees $\rho =1$ and
we are dealing with geometric $p$-rough paths of Brownian regularity, i.e. $%
p=2+\varepsilon $.) Condition (\ref{CondOfStableCorollary}) then holds with $%
H=1/2:$ take $s,t\in \left[ 0,1/n\right] $ and compute, with focus on one
non-diagonal entry,%
\begin{eqnarray*}
R_{X-B}\left( s,t\right)  &\equiv &\mathbb{E[}\left( X_{s}-B_{s}\right)
\left( X_{t}-B_{t}\right) ] \\
&=&\int_{0}^{t}\left( e^{-\left( s-r\right) }-1\right) \left( e^{-\left(
t-r\right) }-1\right) dr=O\left( n^{-3}\right) .
\end{eqnarray*}%
By corollary \ref{FV07cor} we see condition \ref{CondLikeFBMatZero} holds
for the Ornstein Uhlenbeck examples.
\end{example}

\begin{example}[Brownian Bridge]
Writing $X_{t}^{T}:=B_{t}-\frac{t}{T}B_{T}$ where $B$ is standard Brownian
motion this follows along the same lines, again by comparison with $B_{t}$
for $t\rightarrow 0+$.
\end{example}

\section{Taylor Expansions for Rough Differential Equations}

Given a smooth vector field $W$ and smooth driving signal $x\left( \cdot
\right) $ for the ODE $dy=V\left( y\right) dx$, it follows from (\ref%
{LieIdentity}) that%
\begin{equation*}
J_{0\leftarrow t}^{x}\left( W\left( y_{t}^{x}\right) \right) =W\left(
y_{0}\right) +\int_{0}^{t}J_{0\leftarrow s}^{x}\left( \left[ V_{i},W\right]
\left( y_{s}^{x}\right) \right) dx_{s}^{i},
\end{equation*}%
where Einstein's summation convention is used throughout. Iterated use of
this leads to the Taylor expansion%
\begin{eqnarray*}
J_{0\leftarrow t}^{x}\left( W\left( y_{t}^{x}\right) \right)  &=&W|_{y_{0}}+
\left[ V_{i},W\right] |_{y_{0}}\mathbf{x}_{0,t}^{1;i} \\
&&+\left[ V_{i},\left[ V_{j},W\right] \right] |_{y_{0}}\mathbf{x}%
_{0,t}^{2;i,j} \\
&&+\dots  \\
&&+\left[ V_{i_{1}},\dots \left[ V_{i_{N}},W\right] \right] |_{y_{0}}\mathbf{%
x}_{0,t}^{N;i_{1},\dots ,i_{N}} \\
&&+\cdots 
\end{eqnarray*}%
where $\mathbf{x}_{0,t}$ denotes the signature of $x\left( \cdot \right) |_{%
\left[ 0,t\right] }^{i}$ in $\mathbb{R}^{d}\oplus \left( \mathbb{R}%
^{d}\right) ^{\otimes 2}\oplus \dots \oplus \left( \mathbb{R}^{d}\right)
^{\otimes N}\oplus \dots $ Such an expansion makes immediate sense when $x$
is replaced by a geometric $p$-rough path\footnote{%
By definition, such a $p$-rough path takes values in the step-$\left[ p%
\right] $ tensor algebra but recall that there is a unique lift to the step-$%
N$ group for any $N>\left[ p\right] $.}. Remainder estimates can be obtained
via Euler-estimates \cite{FV07c} provided $J_{0\leftarrow t}^{x}\left(
W\left( y_{t}^{x}\right) \right) $ is a solution of some ODE of form $dz=%
\hat{V}\left( z\right) dx$. This is accomplished by setting%
\begin{equation*}
z:=\left( z^{1},z^{2},z^{3}\right) :=\left( y^{x},J_{0\leftarrow
t}^{x},J_{0\leftarrow t}^{x}\left( W\left( y_{t}^{x}\right) \right) \right)
\in \mathbb{R}^{e}\oplus \mathbb{R}^{e\times e}\oplus \mathbb{R}^{e}
\end{equation*}%
Noting that $J_{0\leftarrow t}^{x}\left( W\left( y_{t}^{x}\right) \right) $
is given by $z^{2}\cdot W\left( z^{1}\right) $ in terms of matrix
multiplication we have 
\begin{eqnarray*}
dz^{1} &=&V_{i}\left( z^{1}\right) dx^{i} \\
dz^{2} &=&-z^{2}\cdot DV_{i}\left( z^{1}\right) dx^{i} \\
dz^{3} &=&\left( dz^{2}\right) \cdot W\left( z^{1}\right) +z^{2}\cdot
d\left( W\left( z^{1}\right) \right)  \\
&=&z^{2}\cdot \left( -DV_{i}\left( z^{1}\right) \cdot W\left( z^{1}\right)
+DW\left( z^{1}\right) \cdot V_{i}\left( z^{1}\right) \right) dx^{i} \\
&=&z^{2}\cdot \left[ V_{i},W\right] |_{z^{1}}dx^{i}
\end{eqnarray*}%
started from $\left( y_{0},I,W\left( y_{0}\right) \right) $ where $I$
denotes the identity matrix in $\mathbb{R}^{e\times e}$ and we see that $%
\hat{V}$ is given by%
\begin{equation*}
\hat{V}_{i}\left( z^{1},z^{2},z^{3}\right) =\left( 
\begin{array}{c}
V_{i}\left( z^{1}\right)  \\ 
-z^{2}\cdot DV_{i}\left( z^{1}\right)  \\ 
z^{2}\cdot \left[ V_{i},W\right] \left( z_{1}\right) 
\end{array}%
\right) ,\,\,\,\,i=1,...,d.
\end{equation*}

\begin{lemma}
Assume $V_{1},...,V_{d},W$ are smooth vector fields, bounded with all
derivatives bounded. Then $\hat{V}=\left( \hat{V}_{1},...,\hat{V}_{d}\right) 
$ is a collection of smooth (possibly unbounded) vector fields but explosion
does not occur. More precisley, there exists a unique RDE solution to $dz=%
\hat{V}\left( z\right) d\mathbf{x}$ on any compact time interval $\left[ 0,T%
\right] $. In fact, for some\ increasing function $\varphi $ from $\mathbb{R}%
^{+}$ into itself%
\begin{equation*}
\left\vert z\right\vert _{\infty ;\left[ 0,t\right] }\leq \varphi \left(
M\right) \text{ when }\left\Vert \mathbf{x}\right\Vert _{p\text{-var;}\left[
0,t\right] }\leq M\text{.}
\end{equation*}
\end{lemma}

\begin{proof}
Smoothness of $\hat{V}$ is obvious and so the RDE $dz=\hat{V}\left( z\right)
d\mathbf{x}$ has a solution up to some possible explosion time. From the
particular structure of $\hat{V}$ we now argue that explosion cannot occur
in finite time: $z^{1}$ does not explode as it is a genuine RDE solution
along Lip vector fields, $z^{2}$ does not explode as it satisfies a linear
RDE (driven by some rough path $M^{\mathbf{x}}$ as already remarked in the
proof of Proposition \ref{VarOfConst}). Clearly then, $z^{3}=z^{2}\cdot
W\left( z^{1}\right) $ where $W$ is a bounded vector fields cannot explode.
More precisely, using the estimates for RDE solutions along $Lip$
respectively linear vector fields in \cite{FV07c} respectively \cite{L98} it
is clear that $z$ remains in a ball of radius only depending on $M$ if $%
\left\Vert \mathbf{x}\right\Vert _{p\text{-var;}\left[ 0,t\right] }\leq M$.
(With some care one can show that $\log \varphi \left( M\right) =O(M^{p})$
as $M\rightarrow \infty $ but this is irrelevant for the sequel.)
\end{proof}

Let us make the following definitions: given $\left( m-1\right) $-times
differentiable vector fields $V=\left( V_{1},...,V_{d}\right) $ on $\mathbb{R%
}^{e}$, $\mathbf{g}\in \oplus _{k=0}^{m}\left( \mathbb{R}^{d}\right)
^{\otimes k}$ and $y\in \mathbb{R}^{e}$ we write%
\begin{equation*}
\mathcal{E}_{\left( V\right) }\left( y,\mathbf{g}\right)
:=\sum_{k=1}^{m}\sum _{\substack{ i_{1},...,i_{k}  \\ \in \left\{
1,...,d\right\} }}\mathbf{g}^{k,i_{1},\cdots ,i_{k}}V_{i_{1}}\cdots
V_{i_{k}}I\left( y\right) .
\end{equation*}%
(Here $I$ denote the identitify function on $\mathbb{R}^{e}$ and vector
fields identified with first order differential operators.) In a similar
spirit, given another sufficiently smooth vector field $W$ we first write%
\begin{equation*}
\lbrack V_{i_{1}},V_{i_{2}},\cdots V_{i_{k}},W]:=[V_{i_{1}},\left[
V_{i_{2}},\cdots \left[ V_{i_{k}},W]\cdots \right] \right]
\end{equation*}%
(which may be viewed as first order differential operator) and then%
\begin{equation}
\mathbf{g}^{k}\cdot \left[ V,\dots ,V,W\right] |_{y_{0}}:=\sum_{\substack{ %
i_{1},...,i_{k}  \\ \in \left\{ 1,...,d\right\} }}\mathbf{g}^{k,i_{1},\cdots
,i_{k}}[V_{i_{1}},V_{i_{2}},\cdots V_{i_{k}},W]
\label{defOfContrLieBrackets}
\end{equation}%
with the convention that $\mathbf{g}^{0}\cdot V_{k}=V_{k}$.

\begin{proposition}
(Localized Euler Estimates) Given $\left\Vert \mathbf{x}\right\Vert _{p\text{%
-var;}\left[ 0,t\right] }\leq M$ and some integer $m>p-1$ there exists $%
C=C\left( M\right) =C=C\left( M,m,p,\hat{V}\right) $ such that%
\begin{equation*}
\left\vert \pi _{\hat{V}}\left( z_{0},0;\mathbf{x}\right) _{0,t}-\mathcal{E}%
_{\left( \hat{V}\right) }\left( z_{0},S_{m}\left( \mathbf{x}\right)
_{0,t}\right) \right\vert \leq C\left( M\right) \times t^{\frac{m+1}{p}}.
\end{equation*}
\end{proposition}

\begin{proof}
If $\hat{V}$ were bounded with bounded derivatives this would be a
consequence of \cite[Thm 19]{FV07c}. On the other hand, $z$ must remain in
the ball $B\left( 0,\varphi \left( M\right) \right) $ and we can replace $%
\hat{V}$ by (compactly supported) vector fields $\tilde{V}$ such that $\hat{V%
}\equiv \tilde{V}$ on $B\left( 0,\varphi \left( M\right) \right) $. After
this localization we apply \cite[Thm 19]{FV07c}.
\end{proof}

\begin{lemma}
Let $f$ be a smooth function on $\mathbb{R}^{e}$ lifted to a smooth function
on $\mathbb{R}^{e}\oplus \mathbb{R}^{e\times e}\oplus \mathbb{R}^{e}$ by%
\begin{equation*}
\hat{f}\left( z^{1},z^{2},z^{3}\right) =f\left( z^{3}\right) .
\end{equation*}
Viewing vector fields as first order differential operators, we have%
\begin{equation*}
\hat{V}_{i_{1}}\cdots \hat{V}_{i_{N}}|_{z_{0}}\hat{f}=\left[
V_{i_{1}},\cdots ,V_{i_{N}},W\right] |_{y_{0}}f.
\end{equation*}%
As a consequence, if $I$ denotes the identity function on $\mathbb{R}^{e}$, 
\begin{eqnarray*}
&&\left\vert z_{t}^{3}-W|_{y_{0}}-\sum_{k=1}^{m}\mathbf{X}_{0,t}^{k}\cdot %
\left[ V,\cdots ,V,W\right] |_{y_{0}}\,\right\vert \\
&\leq &\left\vert z_{0,t}-\mathcal{E}_{\left( \hat{V}\right) }\left(
z_{0},S_{m}\left( \mathbf{x}\right) _{0,t}\right) \right\vert .
\end{eqnarray*}
\end{lemma}

\begin{proof}
Taylor expansion of the evolution ODE of $z^{3}\left( t\right) $ shows that $%
\hat{V}_{i_{1}}\cdots \hat{V}_{i_{N}}|_{z_{0}}f=\left[ V_{i_{1}},\cdots
,V_{i_{N}},W\right] |_{y_{0}}f$ .
\end{proof}

\begin{corollary}
\label{WarwickPresentp18}Fix $a\in \mathcal{T}_{y_{0}}\mathbb{R}^{e}\cong 
\mathbb{R}^{e}$ with $\left\vert a\right\vert =1$. Then%
\begin{equation*}
\mathbb{P}\left[ \left\vert a^{T}J_{0\leftarrow t}^{x}\left( W\left(
y_{t}^{x}\right) \right) -\sum_{k=0}^{m}a^{T}\left( \mathbf{X}%
_{0,t}^{k}\cdot \left[ V,\cdots ,V,W\right] |_{y_{0}}\right) \right\vert
_{t=1/n}>\frac{\varepsilon }{2}n^{-mH}\right] \rightarrow 0\text{ with }%
n\rightarrow \infty .
\end{equation*}
\end{corollary}

\begin{proof}
We estimate this probability by 
\begin{eqnarray*}
&&\mathbb{P}\left[ \left\vert a^{T}J_{0\leftarrow t}^{x}\left( W\left(
y_{t}^{x}\right) \right) -\sum_{k=0}^{m}a^{T}\left( \mathbf{X}%
_{0,t}^{k}\cdot \left[ V,\cdots ,V,W\right] |_{y_{0}}\right) \right\vert
_{t=1/n}>\frac{\varepsilon }{2}n^{-mH};\left\Vert \mathbf{x}\right\Vert _{p%
\text{-var;}\left[ 0,1/n\right] }\leq 1\right] \\
&&+\mathbb{P}\left[ \left\Vert \mathbf{x}\right\Vert _{p\text{-var;}\left[
0,1/n\right] }\geq 1\right] \text{ \ \ \ \ \ \ \ \ \ \ \ \ then, using }%
\left\vert a\right\vert =1\text{ and the previous lemma,} \\
&\leq &\mathbb{P}\left[ \left\vert z_{0,1/n}-\mathcal{E}_{\left( \hat{V}%
\right) }\left( z_{0},S_{m}\left( \mathbf{x}\right) _{0,1/n}\right)
\right\vert >\frac{\varepsilon }{2}n^{-mH};\left\Vert \mathbf{x}\right\Vert
_{p\text{-var;}\left[ 0,1/n\right] }\leq 1\right] +o\left( 1\right) \\
&\leq &\mathbb{P}\left[ C\left( 1\right) \times \left( \frac{1}{n}\right) ^{%
\frac{m+1}{p}}>\frac{\varepsilon }{2}n^{-mH}\right] +o\left( 1\right) \text{
using the localized Euler estimates.}
\end{eqnarray*}%
The probability of the (deterministic) event%
\begin{equation*}
C\left( 1\right) \left( \frac{1}{n}\right) ^{\frac{m+1}{p}}>\frac{%
\varepsilon }{2}\left( \frac{1}{n}\right) ^{mH}
\end{equation*}%
will be zero for $n$ large enough provided $\frac{m+1}{p}>mH$ which is the
case since $p\geq 1,H\leq 1$.
\end{proof}

\section{On H\"{o}rmander's condition\label{OnHcond}}

Let $V=\left( V_{1},...,V_{d}\right) $ denote a collection of smooth vector
fields defined in a neighbourhood of $y_{0}\in $ $\mathbb{R}^{e}$. Given a
multi-index $I=\left( i_{1},...,i_{k}\right) \in \left\{ 1,...,d\right\} ^{k}
$, with length $\left\vert I\right\vert =k$, the vector field $V_{I}$ is
defined by iterated Lie brackets%
\begin{equation}
V_{I}:=[V_{i_{1}},V_{i_{2}},...,V_{i_{k}}]\equiv \lbrack
V_{i_{1}},[V_{i_{2}},...,[V_{i_{k-1}},V_{i_{k}}]...].
\end{equation}%
If $W$ is another smooth vector field defined in a neighbourhood of $%
y_{0}\in $ $\mathbb{R}^{e}$ we write\footnote{%
We introduced this notation already in the previous section, cf. (\ref%
{defOfContrLieBrackets}).} 
\begin{equation*}
\underset{\in \left( \mathbb{R}^{d}\right) ^{\otimes \left( k-1\right) }}{%
\underbrace{a}}\cdot \underset{\text{length }k}{\,\underbrace{\left[ V,\dots
,V,W\right] }}\,:=\sum_{\substack{ i_{1},...,i_{k-1} \\ \in \left\{
1,...,d\right\} }}%
a^{i_{1},...,i_{k-1}}[V_{i_{1}},V_{i_{2}},...,V_{i_{k-1}},W]
\end{equation*}%
Recall that the step-$r$ free nilpotent group with $d$ generators, $G^{r}(%
\mathbb{R}^{d})$,  was realized as submanifold of the tensor algebra 
\begin{equation*}
T^{\left( r\right) }\left( \mathbb{R}^{d}\right) \equiv \oplus
_{k=0}^{r}\left( \mathbb{R}^{d}\right) ^{\otimes k}.
\end{equation*}

\begin{definition}
Given $r\in \mathbb{N}$ we say that condition (H)$_{r}$ holds at $y_{0}\in 
\mathbb{R}^{e}$ if%
\begin{equation}
\mathrm{span}\left\{ V_{I}|_{y_{0}}:\left\vert I\right\vert \leq r\right\} =%
\mathcal{T}_{y_{0}}\mathbb{R}^{e}\cong \mathbb{R}^{e};  \label{Hr}
\end{equation}%
Similarly, we say that (HT)$_{r}$ holds at $y_{0}$ if the span of 
\begin{equation}
\left\{ \pi _{k-1}\left( \mathbf{g}\right) \cdot \underset{\text{ length }k}{%
\underbrace{\left[ V,...,V,V_{i}\right] }|_{y_{0}}}:k=1,...,r;i=1,...,d,%
\text{ }\mathbf{g}\in G^{r-1}(\mathbb{R}^{d})\right\} .  \label{HTr}
\end{equation}%
equals $\mathcal{T}_{y_{0}}\mathbb{R}^{e}\cong \mathbb{R}^{e}$. H\"{o}%
rmander's condtion (H) is satisfied at $y_{0}$ iff (H)$_{r}$ holds for some $%
r\in \mathbb{N}$. Similarly, we say that the H\"{o}rmander-type condtion
(HT) is satisfied at $y_{0}$ iff (HT)$_{r}$ holds for some $r\in \mathbb{N}$%
. (When no confusion arises we omit reference to $y_{0}$.)
\end{definition}

\begin{proposition}
For any fixed $r\in \mathbb{N}$, the span of (\ref{Hr}) equals the span of (%
\ref{HTr}). Consequently, H\"{o}rmander's condition (H) at $y_{0}$ is
equivalent to the H\"{o}rmander-type condtion (HT) at $y_{0}$.
\end{proposition}

\begin{proof}
Given a multi-index $I=\left( i_{1},...,i_{k-1},i_{k}\right) $ of length $%
k\leq r$ and writing $e_{1},...e_{d}$ for the canonical basis of $\mathbb{R}%
^{d}$%
\begin{eqnarray*}
\mathbf{g} &=&\mathbf{g}\left( t_{_{1}},\dots ,t_{k-1}\right) \\
&=&\exp \left( t_{_{1}}e_{i_{1}}\right) \otimes \dots \otimes \exp \left(
t_{k-1}e_{i_{k-1}}\right) \\
&\in &G^{r-1}(\mathbb{R}^{d})\subset T^{r-1}\left( \mathbb{R}^{d}\right) 
\text{.}
\end{eqnarray*}%
(Recall that $T^{r-1}\left( \mathbb{R}^{d}\right) $ is a tensor algebra with
multiplication $\otimes $,~$\exp $ is defined by the usual series and the
CBH formula shows that the so-defined $g$ is indeed in $G^{r-1}(\mathbb{R}%
^{d})$ as claimed.) It follows that any%
\begin{equation*}
\pi _{k-1}\left( \mathbf{g}\right) \cdot \underset{\text{ length }k}{%
\underbrace{\left[ V,...,V,V_{i_{k}}\right] }|_{y_{0}}}
\end{equation*}%
lies in the (HT)$_{r}$-span i.e. the linear span of (\ref{HTr}). Now, the
(HT)$_{r}$-span is a closed linear subspace of $\mathcal{T}_{y_{0}}\mathbb{R}%
^{e}\cong \mathbb{R}^{e}$ and so it is clear that any element of form%
\begin{equation*}
\pi _{k-1}\left( \partial _{\alpha }\mathbf{g}\right) \cdot \underset{\text{
length }k}{\underbrace{\left[ V,...,V,V_{i_{k}}\right] }|_{y_{0}}}
\end{equation*}%
where $\partial _{\alpha }$ stands for any higher order partial derivative
with respect to $t_{1},...,t_{k-1}$ i.e.%
\begin{equation*}
\partial _{\alpha }=\left( \frac{\partial }{\partial t_{1}}\right) ^{\alpha
_{1}}\dots \left( \frac{\partial }{\partial t_{k-1}}\right) ^{\alpha _{k-1}}%
\text{ \ \ with }\alpha \in \left( \mathbb{N}\cup \left\{ 0\right\} \right)
^{k-1}
\end{equation*}%
is also in the (HT)$_{r}$-span for any $t_{1},\dots ,t_{k-1}$ and, in
particular, when evaluated at $t_{1}=\cdots =t_{k-1}=0$. For the particular
choice $\alpha =\left( 1,\dots ,1\right) $ we have%
\begin{equation*}
\frac{\partial ^{k-1}}{\partial t_{1}\dots \partial t_{_{k-1}}}\mathbf{g}%
|_{t_{1}=0....t_{k-1}=0}=e_{i_{1}}\otimes \dots \otimes e_{i_{k-1}}=:\mathbf{%
h}
\end{equation*}%
where\ $\mathbf{h}$ is an element of $T^{r-1}\left( \mathbb{R}^{d}\right) $
with the only non-zero entry arising on the $\left( k-1\right) ^{th}$ tensor
level, i.e.%
\begin{equation*}
\pi _{k-1}\left( \mathbf{h}\right) =e_{i_{1}}\otimes \dots \otimes
e_{i_{k-1}}.
\end{equation*}%
Thus,%
\begin{equation*}
\pi _{k-1}\left( \mathbf{h}\right) \cdot \underset{\text{ length }k}{%
\underbrace{\left[ V,...,V,V_{i_{k}}\right] }|_{y_{0}}}=\left[
V_{i_{1}},\dots ,V_{i_{k-1}},V_{i_{k}}\right] |_{y_{0}}
\end{equation*}%
is in our (HT)$_{r}$- span. But this says precisely that, for any
multi-index $I$ of lenght $k\leq r$ the bracket vector field evaluated at $%
y_{0}$ i.e. $V_{I}|_{y_{0}}$ is an element of our (HT)$_{r}$-span.
\end{proof}

\section{Proof of Main Result}

We are now in a position to give

\begin{proof}[Proof (of Theorem \protect\ref{Main})]
We fix $t\in (0,T]$. As usual it suffices to show a.s. invertibility of%
\begin{equation*}
\sigma _{t}=\left( \left\langle DY_{t}^{i},DY_{t}^{j}\right\rangle _{%
\mathcal{H}}\right) _{i,j=1,...,e}\in \mathbb{R}^{e\times e}.
\end{equation*}%
In terms of an ONB $\left( h_{n}\right) $ of the Cameron Martin space we can
write%
\begin{eqnarray}
\sigma _{t} &=&\sum_{n}\left\langle DY_{t},h_{n}\right\rangle _{\mathcal{H}%
}\otimes \left\langle DY_{t},h_{n}\right\rangle _{\mathcal{H}}
\label{MCovMatrix} \\
&=&\sum_{n}\int_{0}^{t}J_{t\leftarrow s}^{\mathbf{X}}\left( V_{k}\left(
Y_{s}\right) \right) dh_{n,s}^{k}\otimes \int_{0}^{t}J_{t\leftarrow s}^{%
\mathbf{X}}\left( V_{l}\left( Y_{s}\right) \right) dh_{n,s}^{l}  \notag
\end{eqnarray}%
(Summation over up-down indices is from here on tacitly assumed.)
Invertibility of $\sigma $ is equivalent to invertibility of the reduced
covariance matrix%
\begin{equation*}
C_{t}:=\sum_{n}\int_{0}^{t}J_{0\leftarrow s}^{\mathbf{X}}\left( V_{k}\left(
Y_{s}\right) \right) dh_{n,s}^{k}\otimes \int_{0}^{t}J_{0\leftarrow s}^{%
\mathbf{X}}\left( V_{l}\left( Y_{s}\right) \right) dh_{n,s}^{l}
\end{equation*}%
which has the advantage of being adapted, i.e. being $\sigma \left(
X_{s}:s\in \left[ 0,t\right] \right) $-measurable. We now assume that%
\begin{equation*}
\mathbb{P}\left( \det C_{t}=0\right) >0
\end{equation*}%
and will see that this leads to a contradiction with H\"{o}rmander's
condition.\newline
\underline{Step 1:} Let $K_{s}$ be the random subspace of $\mathcal{T}%
_{y_{0}}\mathbb{R}^{e}\cong \mathbb{R}^{e}$. spanned by%
\begin{equation*}
\left\{ J_{0\leftarrow r}^{\mathbf{X}}\left( V_{k}\left( Y_{r}\right)
\right) ;r\in \left[ 0,s\right] ,k=1,...,d\right\} .
\end{equation*}%
The subspace $K_{0^{+}}=\cap _{s>0}K_{s}$ is measurable with respect to the
germ $\sigma $-algebra and by our "0-1 law" assumption, deterministic with
probability one. A random time is defined by%
\begin{equation*}
\Theta =\inf \left\{ s\in (0,t]:\dim K_{s}>\dim K_{0^{+}}\right\} \wedge t,
\end{equation*}%
and we note that $\Theta >0$ a.s. For any vector $v\in \mathbb{R}^{e}$ we
have%
\begin{equation*}
v^{T}C_{t}v=\sum_{n}\left\vert \int_{0}^{t}v^{T}J_{0\leftarrow s}^{\mathbf{X}%
}\left( V_{k}\left( Y_{s}\right) \right) dh_{n,s}^{k}\right\vert ^{2}.
\end{equation*}%
Assuming $v^{T}C_{t}v=0$ implies%
\begin{equation*}
\forall n:\int_{0}^{t}v^{T}J_{0\leftarrow s}^{\mathbf{X}}\left( V_{k}\left(
Y_{s}\right) \right) dh_{n,s}^{k}=0
\end{equation*}%
and hence, by our non-degeneracy condition on the Gaussian process%
\begin{equation*}
v^{T}J_{0\leftarrow s}^{\mathbf{X}}\left( V_{k}\left( Y_{s}\right) \right) =0
\end{equation*}%
for any $s\in \lbrack 0,t]$ and any $k=1,...,d$ which implies that $v$ is
orthogonal to $K_{t}$. Therefore, $K_{0^{+}}\neq \mathbb{R}^{e}$, otherwise $%
K_{s}=\mathbb{R}^{e}$ for every $s>0$ so that $v$ must be zero, which
implies $C_{t}$ is invertible a.s. in contradiction with our hypothesis.%
\newline
\underline{Step 2:} We saw that $K_{0^{+}}$ is a deterministic and linear
subspace of $\mathbb{R}^{e}$ with strict inclusion $K_{0^{+}}\subsetneqq 
\mathbb{R}^{e}$ In particular, there exists a deterministic vector $z\in $ $%
\mathbb{R}^{e}\backslash \left\{ 0\right\} $ which is orthogonal to $%
K_{0^{+}}.$ We will show that $z$ is orthogonal to to all vector fields and
(suitable) brackets evaluated at $y_{0}$, thereby contradicting the fact
that our vector fields satisfy H\"{o}rmander's condition. By definition of $%
\Theta $, $K_{0^{+}}\equiv K_{t}$ for $0\leq t<\Theta $ and so for every $%
k=1,...d,$%
\begin{equation}
z^{T}J_{0\leftarrow t}^{\mathbf{X}}\left( V_{k}\left( Y_{t}\right) \right) =0%
\text{ for }t\leq \Theta .  \label{zOrth}
\end{equation}

\noindent Observe that, by evaluation at $t=0$, this implies $z$ $\bot $ span%
$\left\{ V_{1},...,V_{d}\right\} |_{y_{0}}$.

\underline{Step 3:}\noindent\ We call an element $\mathbf{g}\in \oplus
_{k=0}^{\infty }\left( \mathbb{R}^{d}\right) ^{\otimes k}$ group-like iff
for any $N\in \mathbb{N}$,%
\begin{equation*}
\left( \pi _{0}\left( \mathbf{g}\right) ,\dots ,\pi _{N}\left( \mathbf{g}%
\right) \right) \in G^{N}\left( \mathbb{R}^{d}\right) \subset \oplus
_{k=0}^{N}\left( \mathbb{R}^{d}\right) ^{\otimes k}.
\end{equation*}%
We now keep $k$ fixed and make \textbf{induction hypothesis} $I\left(
m-1\right) :$%
\begin{equation*}
\forall \mathbf{g}\text{ group-like, }j\leq m-1:\text{ }z^{T}\pi _{j}\left( 
\mathbf{g}\right) \left. \left[ V,\cdots ,V;V_{k}\right] \right\vert
_{y_{0}}=0.
\end{equation*}%
To this end, take the shortest path $\gamma ^{n}:\left[ 0,1/n\right]
\rightarrow \mathbb{R}^{d}$ such that $S_{m}\left( \gamma ^{n}\right) $
equals $\pi _{1,...,m}\left( \mathbf{g}\right) $, the projection of $\mathbf{%
g}$ to the free step-$m$ nilpotent group with $d$ generators, denoted $%
G^{m}\left( \mathbb{R}^{d}\right) $. Then%
\begin{equation*}
\left\vert \gamma ^{n}\right\vert _{1\text{-var;}\left[ 0,1/n\right]
}=\left\Vert \pi _{1,...,m}\left( \mathbf{g}\right) \right\Vert
_{G^{m}\left( \mathbb{R}^{d}\right) }<\infty 
\end{equation*}%
and the scaled path%
\begin{equation*}
h^{n}\left( t\right) =n^{-H}\gamma ^{n}\left( t\right) ,\,\,\,H\in \left(
0,1\right) 
\end{equation*}%
has length (over the interval $\left[ 0,1/n\right] $) proportional to $n^{-H}
$ which tends to $0$ as $n\rightarrow \infty $. Our plan is to show that%
\begin{equation}
\forall \varepsilon >0:\underset{n\rightarrow \infty }{\lim \inf }\mathbb{P}%
\left( \left\vert z^{T}J_{0\leftarrow 1/n}^{h^{n}}\left( V_{k}\left(
y_{1/n}^{h^{n}}\right) \right) \right\vert <\varepsilon /n^{mH}\right) >0
\label{ToShowInStep3}
\end{equation}%
which, since the event involved is deterministic, really says that%
\begin{equation*}
\left\vert n^{mH}z^{T}J_{0\leftarrow 1/n}^{h^{n}}\left( V_{k}\left(
y_{1/n}^{h^{n}}\right) \right) \right\vert <\varepsilon 
\end{equation*}%
holds true for all $n\geq n_{0}\left( \varepsilon \right) $ large enough.
Then, sending $n\rightarrow \infty $, a Taylor expansion and $I\left(
m-1\right) $ shows that the l.h.s. converges to%
\begin{equation*}
\left\vert z^{T}\underset{=\pi _{m}\left( \mathbf{g}\right) }{\underbrace{%
n^{mH}\pi _{m}\left( S_{m}\left( h^{n}\right) \right) }}\left. \cdot \left[
V,\cdots ,V;V_{k}\right] \right\vert _{_{y_{0}}}\right\vert <\varepsilon 
\text{ }
\end{equation*}%
and since $\varepsilon >0$ is arbitrary we showed $I\left( m\right) $ which
completes the induction step.

\underline{Step 4:}\noindent\ \ The only thing left to show is (\ref%
{ToShowInStep3}), that is, positivity of $\lim \inf $ of 
\begin{eqnarray*}
&&\mathbb{P}\left( \left\vert z^{T}J_{0\leftarrow 1/n}^{h^{n}}\left(
V_{k}\left( y_{1/n}^{h^{n}}\right) \right) \right\vert <\varepsilon
/n^{mH}\right)  \\
&\geq &\mathbb{P}\left( \left\vert z^{T}J_{0\leftarrow \cdot }^{\mathbf{X}%
}\left( V_{k}\left( y_{\cdot }\right) \right) -z^{T}J_{0\leftarrow \cdot
}^{h^{n}}\left( V_{k}\left( y_{\cdot }^{h^{n}}\right) \right) \right\vert
_{\cdot =1/n}<\varepsilon /n^{mH}\right)  \\
&&-\mathbb{P}\left( \Theta \leq 1/n\right) 
\end{eqnarray*}%
and since $\Theta >0$ a.s. it is enough to show that 
\begin{equation*}
\underset{n\rightarrow \infty }{\lim \inf }\mathbb{P}\left( \left\vert
z^{T}J_{0\leftarrow \cdot }^{\mathbf{X}}\left( V_{k}\left( y_{\cdot }\right)
\right) -z^{T}J_{0\leftarrow \cdot }^{h^{n}}\left( V_{k}\left( y_{\cdot
}^{h^{n}}\right) \right) \right\vert _{\cdot =1/n}<\varepsilon
/n^{mH}\right) >0.
\end{equation*}%
Using $I\left( m-1\right) $ + stochastic Taylor expansion (more precisely,
corollary \ref{WarwickPresentp18}) this is equivalent to show positivity of $%
\lim \inf $ of%
\begin{equation*}
\mathbb{P}\left( \left\vert z^{T}\mathbf{X}_{0,\cdot }^{m}\left[ V,\cdots
,V;V_{k}\right] -z^{T}J_{0\leftarrow \cdot }^{h^{n}}\left( V_{k}\left(
y_{\cdot }^{h^{n}}\right) \right) \right\vert _{\cdot =1/n}<\varepsilon
/n^{mH}\right) .
\end{equation*}%
Rewriting things, we need to show positivity of $\lim \inf $ of%
\begin{equation*}
\mathbb{P}(|n^{mH}z^{T}\left[ V,\cdots ,V;V_{k}\right] \mathbf{X}%
_{0,1/n}^{m}-\underset{\rightarrow z^{T}\left[ V,\cdots ,V;V_{k}\right] \pi
_{m}\left( \mathbf{g}\right) }{\underbrace{z^{T}n^{mH}J_{0\leftarrow
1/n}^{h^{n}}\left( V_{k}\left( y_{1/n}^{h^{n}}\right) \right) }}%
|<\varepsilon )
\end{equation*}%
or, equivalently, 
\begin{equation*}
\mathbb{P}\left( \left\vert z^{T}\left. \left[ V,\cdots ,V;V_{k}\right]
\right\vert _{y_{0}}\left( n^{mH}\mathbf{X}_{0,1/n}^{m}-\pi _{m}\left( 
\mathbf{g}\right) \right) \right\vert <\varepsilon \right) 
\end{equation*}%
But this is implied by condition \ref{CondLikeFBMatZero} and so the proof is
finished.
\end{proof}

\textbf{Acknowledgment: }The second author is partially supported by a
Leverhulme Research Fellowship.\newline

\end{document}